# Pulse-and-Glide Driving with Drivability Constraints: A Pontryagin Approach

Carey A. Whitehair, Michelle A. K. Denlinger, and Hosam K. Fathy
Department of Mechanical and Nuclear Engineering
The Pennsylvania State University
E-mail: hkf2@engr.psu.edu

This paper uses Pontryagin methods to analyze pulse-and-glide driving for different nominal vehicle speeds. There is significant literature supporting the fact that pulse-and-glide has the potential to reduce fuel consumption compared to driving at a constant speed, but the benefits at a variety of nominal speeds remain relatively unexplored. Building on the literature, we formulate a speed trajectory optimization problem where the objective is a linearly-scalarized Pareto combination of fuel consumption and average vehicle speed, in addition to a quadratic penalty on jerk. By analyzing this optimization problem using Pontryagin methods, we show that (i) for each nominal speed, there is a critical penalty on jerk below which the optimal solution is pulse-and-glide, (ii) without any penalty on jerk, the optimal pulse-and-glide trajectory switches infinitely fast, and (iii) above a critical nominal velocity, the optimal solution is steady-speed driving rather than pulse-and-glide, regardless of the penalty on jerk.
Topics / Powertrain Control and Energy Management

## 1. INTRODUCTION

This paper uses Pontryagin methods to analyze the degree to which soft penalties on drivability affect pulse-and-glide (PnG) driving. PnG driving is a "dynamic coasting" policy for adjusting the speed of on-road vehicles as a function of time. It involves switching periodically between two different propulsion power levels, where a period of acceleration is followed by a period of deceleration through coasting. It is analogous to "dynamic soaring" in flying birds, which is known to enable longer flights with lower energy consumption.

The main benefit of PnG driving is a potentially significant reduction in vehicle fuel consumption, compared to constant-speed driving. This benefit stems from the non-convexity of internal combustion engine (ICE) fuel consumption as a function of propulsion power. Figure 1 illustrates this argument. The red dot labeled $Q_{CS}$ indicates the constant fuel rate required to travel at a constant speed. The green dots show a PnG policy that alternates between the minimum brake-specific fuel consumption (BSFC) and idling. As shown, the average PnG fuel consumption rate, indicated by the green dot labeled $Q_{PnG}$, is below the fuel rate for constant speed. The extent to which fuel consumption is reduced depends on the precise shape of the ICE's BSFC map, and the average propulsion power demanded.

There is significant existing literature supporting the fact that PnG driving has the potential to reduce vehicle fuel consumption compared to steady-speed driving. Much of this literature builds on Lee's work on the periodic switching between the minimum BSFC point and turning the engine completely off, thereby reducing drive power to zero [1]. Motivated by this work, the literature examines a number of alternative PnG policies, including keeping the engine in idle [2]–[4] and keeping the engine close to the minimum BSFC point [3]–[4]. Furthermore, studies have examined the efficacy of PnG in both single vehicles [1]–[3] and platoons [6], [7]. PnG driving has also been tested in on-road experiments [1]. Theoretical tools have been used to analyze PnG driving such as the Π-test [4] and Pontryagin methods [8], [9].

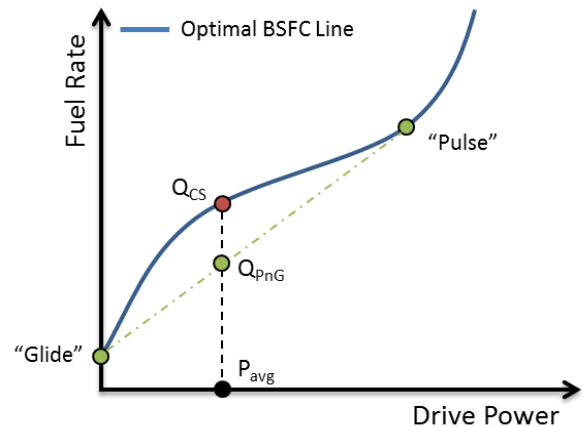

Fig. 1: BSFC as a function of power during PnG.

The above literature, while quite rich, leaves a number of important questions pertaining to PnG driving relatively unexplored. First, the nominal vehicle speed to which PnG driving is compared in the literature is typically the speed that minimizes steady-state fuel consumption. As a result, the benefits of PnG driving around other nominal speeds remain relatively unexplored. Second, PnG driving is known to result in aggressive speed trajectories, which may not be acceptable from the perspective of passenger comfort. This creates a motivation for studying the tradeoffs between fuel consumption minimization and drivability. Our work addresses these two challenges using Pontryagin methods. Specifically, we formulate a vehicle speed trajectory optimization problem where the optimization objective is a linearly-scalarized Pareto combination of fuel consumption minimization and speed maximization. Second, we introduce a soft quadratic penalty on longitudinal jerk. By varying this penalty, we are able to show the existence of a transition point where PnG driving is no longer superior to steady-speed driving.

The remainder of this paper is organized as follows. Section 2 outlines the problem formulation including the vehicle model. Section 3 presents the analysis of the optimal trajectory through the Pontryagin minimum principle. Section 4 presents the linearization of these results around steady state. Last, section 5 validates these results through simulation.

## 2. PROBLEM FORMULATION
### 2.1 Vehicle Model

The formulation of our vehicle speed trajectory optimization considers a 1991 Dodge Caravan as modeled by driving on flat terrain with aerodynamic drag and friction forces opposing the vehicle's motion. In our longitudinal model, we use vehicle velocity and engine propulsive force as states. The state equations are

$$\dot{x}_1 = \frac{1}{M}\left(x_2 - \frac{1}{2}\rho C_d A_f x_1^2 - \mu M g\right) \quad (1)$$
$$\dot{x}_2 = u,$$

where $x_1$ is the vehicle velocity, $x_2$ is the propulsive force, and the input, $u$, is the rate of change of propulsion force with respect to time (i.e., a quantity that represents longitudinal jerk). Constant parameters $M$, $\rho$, $C_d$, $A_f$, $\mu$, and $g$ represent vehicle mass, air density, vehicle drag coefficient, frontal area, coefficient of friction, and gravitational constant, respectively, and are found in Table 1. The vehicle model parameters are obtained from the Advanced Vehicle Simulator (ADVISOR) [10].

Engine power $P$ is modeled as the product of propulsive force and vehicle speed:
$$P = x_1 x_2. \quad (2)$$
BSFC is approximated as a quadratic function of engine power:
$$\beta(P) = \beta_0 + \frac{\gamma}{2}(P - P_0)^2, \quad (3)$$
where constant parameters $\beta_0$, $\gamma$, and $P_0$ are fitted from the data obtained from the ADVISOR vehicle model using the MATLAB curve fitting function 'fit'. The rate of fuel consumption $\dot{m}_f$ is modeled as the product of BSFC and engine power:
$$\dot{m}_f = P\,\beta(P). \quad (4)$$
From (2)-(4) we obtain an approximation of the optimal BSFC curve which minimizes fuel consumption for a given power demand. The vehicle parameters used for this work are shown in Table 1. The fit of the BSFC curve to the raw data is shown in Figure 2. The curve fit assumes zero fuel consumption at idle. In addition, this model assumes an ideal, continuously variable transmission with negligible driveline energy losses. This type of vehicle model is consistent with the existing PnG literature [4],[7].

Table 1 Vehicle Parameters

| Mass | $M$ | 1605 | kg |
|---|---|---|---|
| Air density | $\rho$ | 120 | kg/m$^3$ |
| Frontal Area | $A_f$ | 2 | m$^2$ |
| Coefficient of drag | $C_d$ | 0.33 | |
| Coefficient of friction | $\mu$ | 0.009 | |
| Bsfc curve parameter | $B_0$ | 6.5*10$^{-5}$ | |
| Bsfc curve parameter | $P_0$ | 30000 | |
| Bsfc curve parameter | $\gamma$ | 1.1*10$^{-13}$ | |

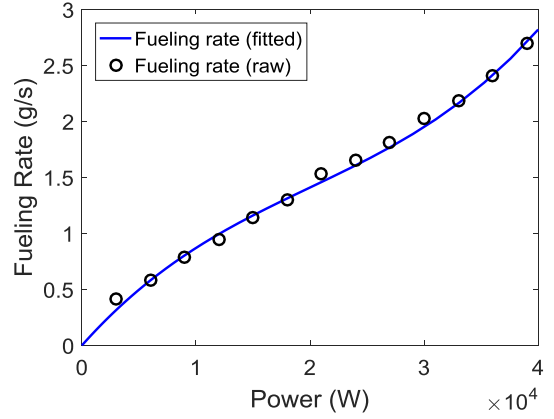

Fig. 2: Optimal BSFC curve

### 2.2 Optimal Control Problem

The optimization objective linearly weights fuel consumption and average vehicle speed, in addition to a quadratic penalty on jerk. In full, the speed trajectory optimization problem is stated as follows:

$$\min_{u,x,T} J(u,x,T) = \frac{1}{T}\int_0^T \dot{m}_f(x) - Cx_1 + \frac{1}{2}Ru^2\,dt \quad (5)$$

s.t. $\dot{m}_f(x) = x_1 x_2 \left(\beta_0 + \frac{\gamma}{2}(x_1 x_2 - P_0)^2\right) \quad (6)$

$$\dot{x}_1 = \frac{1}{M}\left(x_2 - \frac{1}{2}\rho C_d A_f x_1^2 - \mu M g\right) \quad (7)$$
$$\dot{x}_2 = u$$

$$x_2 \geq 0 \quad (8)$$
$$x_1(0) = x_1(T) \quad (9)$$
$$x_2(0) = x_2(T).$$

## 3. APPLICATION OF THE PONTRYAGIN MINIMUM PRINCIPLE

The above optimal control problem can be analyzed using the Pontryagin Minimum Principle (PMP). The Hamiltonian for this problem is given by:

$$H(x,u,\lambda,t) = \dot{m}_f(x) - Cx_1 + \frac{1}{2}Ru^2$$
$$+ \lambda_1 \frac{1}{M}\left(x_2 - \frac{1}{2}\rho C_d A_f x_1^2 - \mu M g\right) + \lambda_2 u \quad (10)$$

where each $\lambda$ is a time-varying Lagrange multiplier or co-state. The dimension of $\lambda$ corresponds to the number of states. The co-state dynamics are given by

$$\dot{\lambda}_1 = -\frac{\partial \dot{m}_f}{\partial x_1} + C + \frac{1}{M}\rho C_d A_f x_1 \lambda_1$$
$$\dot{\lambda}_2 = -\frac{\partial \dot{m}_f}{\partial x_2} - \frac{\lambda_1}{M}. \quad (11)$$

The optimal input is given by

$$\min_u H(x,u,\lambda,t), \quad (12)$$

where, because of the convexity of the quadratic jerk penalty, the optimal input is

$$u = -\frac{\lambda_2}{R}. \quad (13)$$

Necessary boundary conditions for periodic processes like this application are given in [12]. These boundary conditions consist of periodicity constraints and a transversality condition, both defined for a period of $T$. For this problem, these boundary conditions are

$$\lambda_1(0) = \lambda_1(T)$$
$$\lambda_2(0) = \lambda_2(T) \quad (14)$$
$$H(x,u,\lambda,T) - J(u,x,T) = 0 \quad (15)$$

In full, the PMP dynamics and boundary conditions are given by

$$\dot{x}_1 = \frac{1}{M}(x_2 - \frac{1}{2}\rho C_d A_f x_1^2 - \mu M g$$
$$\dot{x}_2 = -\frac{\lambda_2}{R}$$
$$\dot{\lambda}_1 = -\frac{\partial \dot{m}_f}{\partial x_1} + C + \frac{1}{M}\rho C_d A_f x_1 \lambda_1$$
$$\dot{\lambda}_2 = -\frac{\partial \dot{m}_f}{\partial x_2} - \frac{\lambda_1}{M} \quad (16)$$

$$x_1(0) = x_1(T)$$
$$x_2(0) = x_2(T)$$
$$\lambda_1(0) = \lambda_1(T)$$
$$\lambda_2(0) = \lambda_2(T)$$
$$H(x,u,\lambda,T) - J(u,x,T) = 0$$

## 4. LINEARIZATION AND LINEAR ANALYSIS

In this section, we analyze the nature of the optimal solution linearized around an equilibrium point. Specifically, our goal is to analyze the PMP dynamics around different nominal vehicle speeds. To do this, we pick a nominal vehicle speed $x_{1,0}$ and then work backwards to determine the optimization weight $C$ that gives this speed.

For a given nominal speed $x_{1,0}$ the corresponding equilibrium point of the PMP conditions is given by

$$x_{2,0} = \frac{1}{2}\rho C_d A_f x_{1,0}^2 + \mu M g$$
$$\lambda_{1,0} = -M\frac{\partial \dot{m}_f}{\partial x_2} \quad (17)$$
$$\lambda_{2,0} = 0$$

where $\partial \dot{m}_f / \partial x_2$ is evaluated at equilibrium. The corresponding optimization weight $C$ is given by

$$C = \frac{\partial \dot{m}_f}{\partial x_1} + \frac{\partial \dot{m}_f}{\partial x_2}\rho C_d A_f x_{1,0} \quad (18)$$

where $\partial \dot{m}_f/\partial x_1$ and $\partial \dot{m}_f/\partial x_2$ are evaluated at equilibrium. Note that the optimization weight $R$ does not affect the location of the equilibrium.

Linearization is completed around this equilibrium point and weight $C$. This leads to the following system matrix for the linearized PMP conditions.

$$A = \begin{bmatrix} -\frac{\rho C_d A_f x_1}{M} & \frac{1}{M} & 0 & 0 \\ 0 & 0 & 0 & -\frac{1}{R} \\ -\frac{\partial^2 \dot{m}_f}{\partial x_1^2} + \frac{\lambda_1 \rho C_d A_f}{M} & -\frac{\partial^2 \dot{m}_f}{\partial x_1 \partial x_2} & \frac{\rho C_d A_f x_1}{M} & 0 \\ -\frac{\partial^2 \dot{m}_f}{\partial x_1 \partial x_2} & -\frac{\partial^2 \dot{m}_f}{\partial x_2^2} & -\frac{1}{M} & 0 \end{bmatrix}$$
(19)

where the elements in the above $A$ matrix are evaluated at equilibrium.

From (19), we obtain the characteristic equation for the linearized dynamics.

$$0 = s^4 + s^2 \left(\frac{1}{M^2}(\rho C_d A_f x_1)^2 - \frac{1}{R}\frac{\partial^2 \dot{m}_f}{\partial x_2^2}\right)$$
$$+ \frac{1}{R}\frac{1}{M^2}\left(\frac{\partial^2 \dot{m}_f}{\partial x_1^2} + (\rho C_d A_f x_1)^2 \frac{\partial^2 \dot{m}_f}{\partial x_2^2} + \frac{\partial \dot{m}_f}{\partial x_2}\rho C_d A_f\right)$$
(20)

This characteristic equation is dependent on the system parameters, the nominal velocity, and the optimization weight $R$. The roots of this characteristic equation represent the eigenvalues of the linearized PMP conditions. To visualize the impact of the optimization weight R on these roots, we plot a parametric root locus with respect to $R$. For example, the root locus for a nominal speed of 25 m/s is shown in Figure 3.

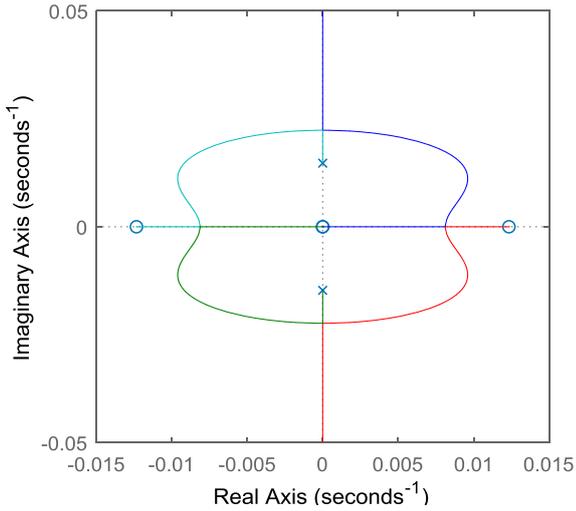

Fig. 3 Root locus at velocity 25 m/s for $R \in [0, \infty)$

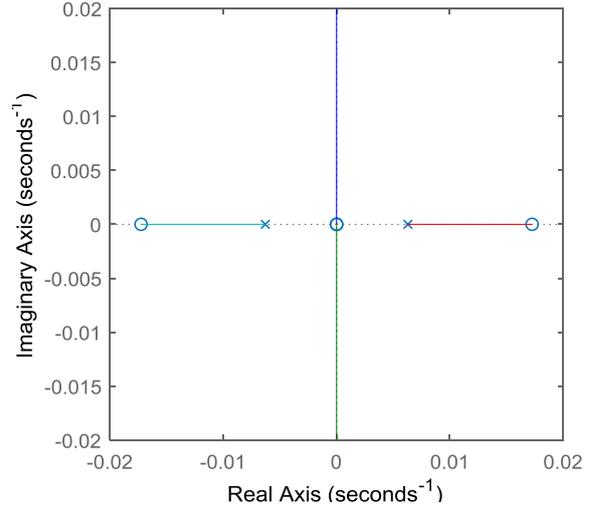

Fig. 4 Root locus at velocity 35 m/s for $R \in [0, \infty)$

The following properties are evident from Figure 3:
- For high values of $R$ (*i.e.*, large penalty on jerk), two of the system's four eigenvalues are unstable. Locally, these unstable eigenvalues are guaranteed to produce trajectories that will *not* satisfy the PMP periodicity conditions. Thus, for large penalty on jerk, the optimal speed trajectory is steady-speed, not PnG.
- For low values of $R$ (*i.e.*, small penalty on jerk), all four of the system's eigenvalues lie on the imaginary axis. For these low $R$ values, the linearized dynamics act like an oscillator, indicating that PnG driving is locally optimal. The frequency of PnG driving is given by the root locus.
- For each nominal driving speed, there is a critical penalty on jerk $R^{\text{crit}}$ for which PnG driving becomes optimal.
- As $R \to 0$, one of the oscillatory frequencies approaches a constant value. The other oscillatory frequency approaches infinity. This result matches what has already been discussed in the PnG literature [11]: when there is no penalty on jerk, the optimal driving pattern is infinitely fast switching.

In addition, we analyze the root locus at several different nominal velocities. Through this analysis, we gain the following insight. *Beyond a critical nominal velocity $v^{\text{crit}}$, the optimal solution will not be PnG, regardless of R value.* For example, Figure 4 shows the root locus for a nominal speed of 35 m/s.

As shown in Figure 4, there is always an unstable eigenvalue, regardless of $R$. For our system, PnG cannot outperform steady state driving for nominal speeds above 33.8 m/s.

Figures 5 and 6 summarize the root locus results for nominal speeds between 2 m/s and 32 m/s. Figure 5 shows how $R^{\text{crit}}$ changes with nominal velocity. In agreement with the findings of [11], the penalty on jerk is relatively small. Figure 6 shows the optimal PnG period at $R^{\text{crit}}$ for different nominal velocities.

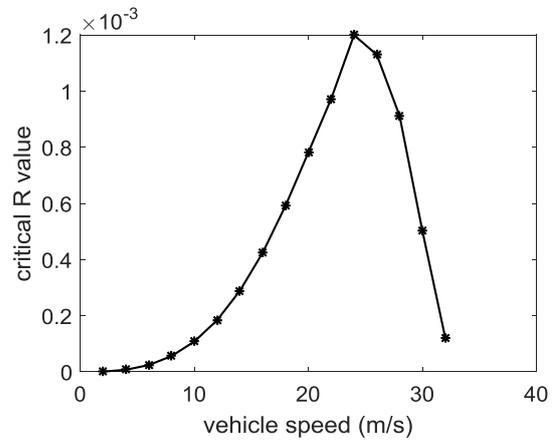

Fig. 5 Critical Values of R

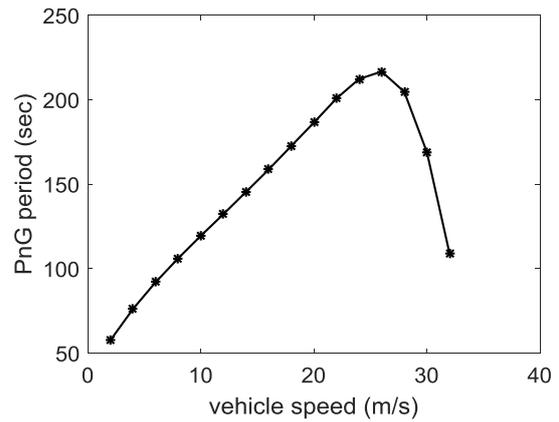

Fig. 6 Period of PnG at $R^{\text{crit}}$

## 5. SIMULATION AND VALIDATION

To validate the results of the linearization, we analyze the nonlinear system. Specifically, we discretize the input and use MATLAB's built-in optimization function 'fmincon' to optimize the initial conditions $x_{1,0}$, $x_{2,0}$ and parameters of the input. The goal is to minimize the cost function and satisfy the constraints listed in (5)-(9). First, we obtain the rough shape of a locally optimal trajectory by parameterizing the input as a single sinusoid. Then, we use those results as an initial condition to obtain a more precise trajectory by parameterizing the input as a Fourier series with six harmonic frequencies.

For both cases, we use cost function weights of $C = 0.42$ (corresponding to a nominal velocity of 15 m/s) and $R = 0.0003$. With these $C$ and $R$ values, the linear analysis predicts oscillatory (PnG) dynamics with two frequencies: $\omega_1 = 0.031$ rad/s and $\omega_2 = 0.0619$ rad/s. Note that the steady-speed driving case results in a cost function of $J_{ss} = -0.2324$.

First, we parameterize the input as a single sinusoid:

$$u = a_1 \sin \omega t + b_1 \cos \omega t, \quad (21)$$

where $\omega$, $a_1$, and $b_1$ are the optimization variables, in addition to initial conditions $x_{1,0}$ and $x_{2,0}$. As an initial guess, we use initial conditions $x_{1,0} = 14.7$ and $x_{2,0}$ equal to the equilibrium value predicted by the linear analysis. The initial guesses for the input are $\omega = 0.1$, $a_1 = 0$, and $b_1 = 7.1$. These initial guesses come close to meeting the periodicity constraints.

The final parameters of the local search performed by 'fmincon' are listed in Table 2, and the optimal trajectory is shown in Figure 7. The cost function is $J = -0.2810$, which outperforms the steady-speed case. As shown in Figure 7, the optimal trajectories deviate significantly from equilibrium. In addition, the optimal solution hits the boundary condition of $x_2 \geq 0$ for the extrema of the simulation. Despite these differences from the linear analysis, the optimal frequency converges to a value very close to the prediction of the linear analysis: 0.0335 rad/s.

Table 2 Optimal Parameters for Single Sinusoid

| | | |
|---|---|---|
| $x_{1,0}$ | 11.97 | m/s |
| $x_{2,0}$ | 230.88 | N |
| $\omega$ | 0.0335 | rad/s |
| $a_1$ | 0.7751 | N/s |
| $b_1$ | 8.4758 | N/s |
| $J$ | $-0.2810$ | |

Next, to gain more insight into the precise shape of the optimal trajectory, we use a 6-harmonic Fourier series to parameterize the input:

$$u = \sum_{k=1}^{6} a_k \sin k\omega t + b_k \cos k\omega t, \quad (22)$$

where $\omega$, $a_1, \ldots, a_6$, $b_1, \ldots, b_6$ are the optimization variables, in addition to initial conditions $x_{1,0}$ and $x_{2,0}$. We perform the optimization several times, first using the single sinusoid trajectory for an initial guess of a 2-harmonic Fourier series and then adding terms until we reach 6 harmonics. The optimal trajectory is shown in Figure 8. This trajectory produces a cost function of $J = -0.3684$, which again outperforms steady-speed driving.

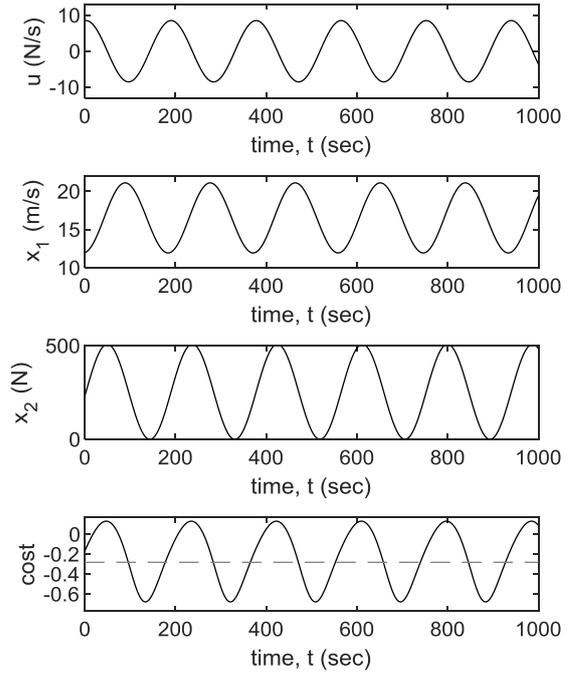

Fig. 7 Optimal trajectory for single sinusoid

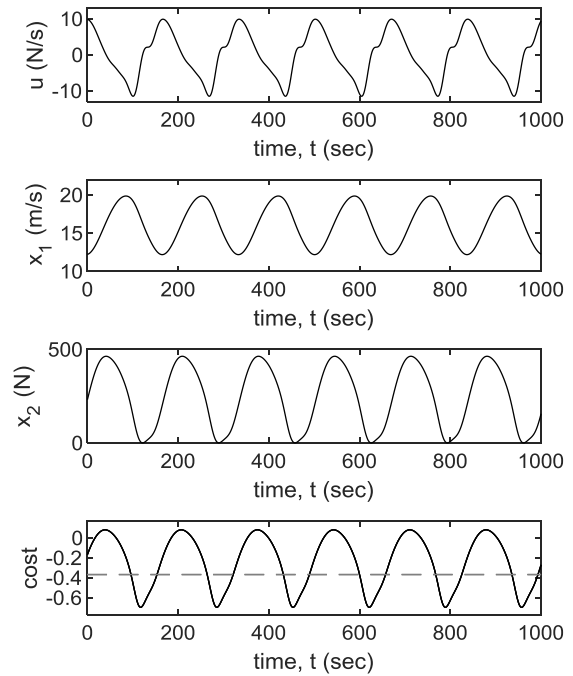

Fig. 8 Optimal trajectory for 6-harmonic Fourier series

## 6. CONCLUSIONS

The analysis in this paper uses Pontryagin methods to determine optimality conditions for the optimization problem listed in (5)–(9). Based on the augmented state and co-state equations, we show the following.

First, for each nominal velocity, there is a critical value of weighting $R$ on jerk below which PnG

outperforms constant speed driving. We show this by examining the root locus of the closed-loop optimal state and co-state equations, linearized about a set of different nominal velocity points. The root locus plot indicates a critical $R$ value, below which the closed-loop poles are purely imaginary. Second, with no penalty on jerk, the optimal solution is PnG with infinitely fast switching. This insight has already been acknowledged by the literature [11]. We validate it using the Pontryagin optimality conditions.

**ACKNOWLEDGEMENTS**

This work was funded in part by the ARPA-E NEXTCAR program and the National Science Foundation CMMI award #1538369, "Self-Adjusting Periodic Optimal Control with Applications to Energy-Harvesting Flight." The authors gratefully acknowledge this support.